\documentclass[11pt]{article}

\usepackage{amssymb, amsmath, amsthm, graphicx}
\usepackage[left=1in,top=1in,right=1in]{geometry}

\usepackage{subfigure}

      \theoremstyle{plain}
      \newtheorem{theorem}{Theorem}

\newtheorem{conjecture}[theorem]{Conjecture}
      \theoremstyle{definition}
      \newtheorem{definition}[theorem]{Definition}

      \theoremstyle{remark}

\def\twr{\mbox{\rm twr}}

\title{Variants of the Erd\H os-Szekeres and Erd\H os-Hajnal Ramsey problems}
\author{Dhruv Mubayi\thanks{Department of Mathematics, Statistics, and Computer Science, University of Illinois, Chicago, IL, 60607 USA.  Research partially supported by NSF grant DMS-1300138. Email: {\tt mubayi@uic.edu}} } 
\begin{document}
\maketitle
\medskip

\begin{abstract}

Given integers $\ell,n$, the $\ell$th power of the  path $P_n$ is the ordered graph $P_n^{\ell}$  with vertex set 
$v_1<v_2<\cdots < v_n$ and all edges of the form $v_iv_j$ where $|i-j|\le \ell$. The  Ramsey number 
$r(P_n^{\ell}, P_n^{\ell})$ is the minimum $N$ such that every 2-coloring of ${[N] \choose 2}$ results in a monochromatic copy of $P_n^{\ell}$.
It is well-known that that $r(P_n^1, P_n^1)=(n-1)^2+1$.   
For $\ell>1$, Balko-Cibulka-Kr\'al-Kyn\v cl proved that $r(P_n^{\ell}, P_n^{\ell})< c_{\ell}n^{128 \ell}$
and asked for the growth rate for fixed $\ell$.
When $\ell=2$, we  improve this upper bound substantially  by proving $r(P_n^{2}, P_n^{2})< cn^{19.5}$.  Using this result, we determine the correct tower growth rate of  the $k$-uniform hypergraph Ramsey number of a $(k+1)$-clique versus an ordered tight path.  Finally, we consider an ordered version of the classical Erd\H os-Hajnal hypergraph Ramsey problem, improve the tower height given by the trivial upper bound,   and conjecture that this tower height is optimal.

\end{abstract}

\section{Introduction}
Let $K_n$ be the complete graph on $n$ vertices.
An ordered path $P_s$ is the graph whose vertices are ordered as $v_1<\cdots  < v_s$ and its edges are $v_1v_2, v_2v_3, \ldots, v_{s-1}v_s$. The Ramsey number $r(P_s, P_n)$ is the minimum $N$ such that every red/blue coloring of ${[N] \choose 2}$ results in a red copy of $P_s$ or a blue copy of $P_n$.
 Note that $P_s$ usually stands for the (unordered) path but since we only consider ordered paths in this paper, we have taken the liberty to use $P_s$ for the ordered path (similarly, $K_n$ refers to the ordered or unordered clique). A similar comment applies to the notation $r(P_s, P_n)$ that we have employed here. Usually this stands for unordered Ramsey numbers but since we exclusively consider ordered Ramsey numbers in this paper, we 
 have chosen to keep this notation and hope that there will  be no confusion. Further, there is not yet a standard notation for ordered Ramsey numbers in the literature as these numbers have only recently  been considered systematically (\cite{BCKK} uses ``${\overline R}$", \cite{CFLS} uses ``$r_{<}$", \cite{MSW} uses ``OR", \cite{MSh} uses ``$N$" and \cite{MS} and the current paper use ``$r$").

Let $f(s,n)$ be the minimum $N$ such that every sequence of $N$ distinct real numbers contains an increasing subsequence of length $s$ or a decreasing subsequence of length $n$. The famous Erd\H os-Szekeres monotone subsequence theorem~\cite{ES35} states that $f(s,n)=(s-1)(n-1)+1$. This function is closely related to $r(P_s, P_n)$. Indeed, it is trivial that $r(P_s, P_n) \ge f(s,n)$ while several proofs for $f(s,n)\le (s-1)(n-1)$ also give the same upper bound for $r(P_s, P_n)$ (see~\cite{MSh} for a further discussion about this). Consequently, it is well-known that 
$r(P_s, P_n)=(s-1)(n-1)+1$. This implies that for fixed $s$, the Ramsey number $r(P_s, P_n)$ is a polynomial function in $n$ (in fact a linear function). On the other hand, $r(n,n)=r(K_n, K_n)$ has exponential growth rate. We begin by considering  
the case of ordered graphs that are denser than  paths but  sparser than cliques.  

\begin{definition} Given $\ell \ge 1$, the $\ell$th power $P_s^{\ell}$ of a path $P_s$ has ordered vertex set $v_1<\cdots  < v_s$ and edge set $\{v_iv_j: |i-j|\le \ell\}$.   The Ramsey number $r(P_s^{\ell}, P_n^{\ell})$ is the minimum $N$ such that every red/blue coloring of ${[N] \choose 2}$ results in a red copy of $P_s^{\ell}$ or a blue copy of $P_n^{\ell}$.
\end{definition}

Note that $P_s^1=P_s$. Conlon-Fox-Lee-Sudakov~\cite{CFLS} asked whether $r(P_n^{\ell}, P_n^{\ell})$ is polynomial in $n$ for every fixed $\ell\ge 1$.  Actually, the problem in~\cite{CFLS} is about the Ramsey number of ordered graphs with bandwidth at most $\ell$ but $P^{\ell}_n$ contains all such graphs so an upper bound for $P^{\ell}_n$ provides an upper bound for the bandwidth problem. This question was answered by
Balko-Cibulka-Kr\'al-Kyn\v cl~\cite{BCKK} who proved an upper bound $c_{\ell}n^{128 \ell}$ and asked~(Problem 2 of~\cite{BCKK}) for the growth rate of $r(P_n^{\ell}, P_n^{\ell})$ (subsequently, a different proof was also sketched in~\cite{CFLS}). The proof of our first result gives a slightly worse polynomial growth rate of $r(P_n^{\ell}, P_n^{\ell})$ than 
in~\cite{BCKK} for large $\ell$ but a much better one for small $\ell$, in particular for $\ell=2$.   Note that the bound in~\cite{BCKK} for $\ell=2$ is $cn^{256}$.
\medskip

\begin{theorem} \label{main}
There is an absolute constant $c>0$ such that
 $r(P_n^{2}, P_n^{2}) < c\, n^{19.487}$ for all $n>1$. 
\end{theorem}

Our second result is an application of Theorem~\ref{main} to a hypergraph Ramsey problem. Indeed, this hypergraph Ramsey problem is what motivated us to consider proving Theorem~\ref{main}.
 
 \begin{definition} A $k$-uniform \emph{tight path} of size $s$, denoted by $P_{s}$, comprises a set of $s$ vertices that are ordered as
 $v_1< \cdots< v_{s}$, and  edges  $(v_j,v_{j+1},\ldots, v_{j + k -1})$ for $j = 1,2,\ldots, s - k + 1$.  The \emph{length} of $P_s$ is the number of edges, $s - k + 1$. Given (ordered) $k$-graphs $F_1, F_2$, the Ramsey number $r_k(F_1,F_2)$ or $r(F_1, F_2)$ is the minimum $N$ such that every red/blue coloring of the edges of the complete $N$-vertex $k$-graph $K^k_N$, whose vertex set is $[N]$, contains a red copy of $F_1$ or a blue copy of $F_2$. 
 \end{definition}

The famous cups-caps theorem of Erd\H os and Szekeres~\cite{ES35}
implies that  $r_3(P_s,P_n) = {n+s - 4 \choose s-2} + 1$. The author and Suk~\cite{MS}
considered the closely related problem of determining $r_k(P_s,n):=r_k(P_s, K_n^k)$ and showed that determining the tower height of $r_4(P_5, n)$ is equivalent to the notorious conjecture of Erd\H os-Hajnal and Rado on the tower height of $r_3(n,n)$. The results in~\cite{MS} focused on  fixed $s$ and large $n$ and there are no nontrivial results for the opposite case, namely for $r_k(s, P_n)$.   Given the close connection between this problem and the problem of determining classical Ramsey numbers, it would be of interest to obtain the growth rate in this range as well. Here we settle the first open case. Recall that the tower function $\twr_i(x)$ is defined by $\twr_1(x)=x$ and $\twr_{i+1}(x)=2^{\twr_i(x)}$. 

\begin{theorem} \label{main2}
For $n$ large, $r_3(4,P_n)< n^{21}$ and more generally for each $k \ge 3$, there  exists  $c>0$ such that for $n$ large, $$\twr_{k-2}(n^{c}) < r_{k}(k+1, P_n) < \twr_{k-2}(n^{62}).$$  
\end{theorem}

The main open problem here is to prove that $r_3(5, P_n)$ has polynomial growth rate and more generally that $r_3(s, P_n)$ has polynomial growth rate for all fixed $s\ge 4$. The corresponding results for higher uniformity follow easily from the case $k=3$.

Our final topic considers a version of the well-known Erd\H os-Hajnal hypergraph Ramsey problem with respect to tight paths. In order to shed more light on classical hypergraph Ramsey numbers,  Erd\H os and Hajnal~\cite{EH72} in 1972 considered the  following more general parameter.
\medskip

\begin{definition} {\bf (Erd\H os-Hajnal~\cite{EH72})} For integers $2\le k < s <n$ and $2 \le t \le {s \choose k}$, let $r_k(s,t;n)$ be the minimum $N$ such that every red/blue coloring of  the edges of $K^k_N$ results in a monochromatic blue copy of $K_n^k$ or has a set of $s$ vertices which induces at least $t$ red edges.
\end{definition}
Note that by definition $r_k(s,n) = r_k(s, {s \choose k}; n)$ so $r_k(s,t;n)$ includes classical Ramsey numbers.
 The main conjecture of Erd\H os and Hajnal states that as $t$ grows from $1$ to ${s\choose k}$, there is a well-defined value $t_1=h_1^{(k)}(s)$ at which $r_k(s,t_1-1;n)$ is polynomial in $n$ while $r_k(s,t_1;n)$ is exponential in a power of $n$, another well-defined value $t_2=h_2^{(k)}(s)$ at which it changes from exponential to double exponential in a power of $n$ and so on,  and finally a well-defined value $t_{k-2}=h_{k-2}^{(k)}(s)<{s \choose k}$
  at which it changes from $\twr_{k-2}$ to $\twr_{k-1}$ in a power of $n$. They were not able to offer a conjecture as to what $h_i^{(k)}(s)$ is in general, except when $i=1$ (for which Erd\H os offered \$500) and when $s=k+1$. For the latter, they conjectured that $h_i^{(k)}(k+1)=i+2$. This was solved for all but three $i$ recently by the author and Suk~\cite{MSEH} via the following result. 

\begin{theorem} {\bf(Mubayi-Suk~\cite{MSEH})} \label{MSEH}
For $4 \le t \le k-2$, there are positive $c =c(k,t)$ and $c'=c'(k,t)$ such that
$$\twr_{t-1}(c' n^{k-t + 1}\log n) \, \ge \, r_k(k+1,t; \, n)
 \, \ge \,  \begin{cases}
\twr_{t-1}(c \, n^{k-t + 1}) \qquad  \hbox{ if $k-t$ is even}\\
\twr_{t-1}(c \, n^{(k-t + 1)/2}) \hskip8pt \hbox{ if $k-t$ is odd.}
\end{cases}
$$
\end{theorem}

 Here we consider the very same problem in the ordered setting by replacing $K_n^k$ with $P_n$.

\begin{definition}   For integers $2\le k < s <n$ and $2 \le t \le {s \choose k}$, let $r_k(s,t;P_n)$ be the minimum $N$ such that every red/blue coloring of  the $k$-sets of $[N]$ results in a monochromatic blue copy of $P_n$ or has a set of $s$ vertices which induces at least $t$ red edges.
\end{definition}

Of course, $r_k(s, {s \choose k}; P_n)= r_k(s, P_n)$.  We will focus our attention on the smallest case $s=k+1$. Our main contribution here is the following conjecture which parallels the Erd\H os-Hajnal conjecture for cliques.

\begin{conjecture} \label{EHpaths}
For $3 \le t \le k$, there are positive $c =c(k,t)$ and $c'=c'(k,t)$ such that
$$\twr_{t-2}(n^c) < r_k(k+1, t; P_n) < \twr_{t-2}(n^{c'}).$$
\end{conjecture}

This conjecture seems more difficult than the original problem of 
Erd\H os and Hajnal. For over 40 years the gaps in the bounds for the 
Erd\H os-Hajnal problem were between exponential and tower functions. Theorem~\ref{MSEH} shows that the correct growth rate is a tower function. For Conjecture~\ref{EHpaths} the gap is again between an exponential  and a tower function but unfortunately the constructions used for Theorem~\ref{MSEH} fail. 

Using standard arguments, it is easy to prove an upper bound of the form $\twr_{t-1}(n^{c})$ in Conjecture~\ref{EHpaths} (see~\cite{MSEH} for details). We improve this upper bound to the tower height given by Conjecture~\ref{EHpaths}.
\bigskip

\begin{theorem}\label{t=k}
For all $3 \le t \le k$ there exists $c=c(k)$ such that $r_k(k+1, t; P_n) < \twr_{t-2}(cn^{2k})$.
\end{theorem}

We make some further modest progress towards Conjecture~\ref{EHpaths} in the cases $t=3$ and $(k,t)=(4,4)$  by proving sharper bounds. Note that Theorem~\ref{main2} determines the correct tower height of $r_k(k+1, t; P_n)$ when $t=k+1$; it is $k-2=t-3$ so the formula differs from that in Conjecture~\ref{EHpaths}. 

\begin{theorem} \label{eh} 
The following bounds hold:

\text{\rm (a)} $r_k(k+1, 3; P_n) \le 16n$ for   $n \ge k \ge 3$

{\rm (b)} $2n-2 \le r_3(4,3; P_n)  \le 3n-4$ for $n \ge 3$

{\rm (c)} $r_4(5, 4; P_n) \ge  2^{n-2}+1$ for  $n\ge 2$.
\end{theorem}

We are not ready to offer a conjecture about the tower growth rate of
$r_k(s, t; P_n)$ as $t$ grows from $2$ to ${s \choose k}$ for $s>k+1$. 

 In all our results where we  find either a long (ordered) blue  $P_n$ or  a small red structure (Theorems\,~\ref{main2},\,\ref{t=k},\,\ref{eh}) our proof actually finds an ordered  blue hypergraph that contains $P_n$. We call this hypergraph a broom (see Definition~\ref{broomdef}). Using brooms we can load the induction hypothesis suitably to carry out the induction step. Perhaps this is one of the main new ideas in this work.

\section{Proof of Theorem~\ref{main}}
We will prove Theorem~\ref{main} by using bounds on the Ramsey multiplicity problem introduced by Erd\H os~\cite{E62}. This problem asks for the largest $\alpha=\alpha(\ell)$ such that every 2-edge coloring of $K_N$
yields at least $(\alpha-o(1)) N^{\ell}$ monochromatic copies of $K_{\ell}$. Erd\H os~\cite{E62} observed 
that $\alpha({\ell})>0$ for all $\ell \ge 4$. The best known bounds for $\alpha({\ell})$ for large $\ell$ are very far apart and can  be found in~\cite{C, T}. We will use the specific result $\alpha(4) > 0.0287/4!=0.0011958\overline{3}$ that was recently proved using Flag Algebras in~\cite{N} (see also \cite{S}). Note that $\log_2(1/\alpha(4))<9.7434$. After this paper was written, we learned that the approach of using ramsey multiplicity for  ordered Ramsey problems was also used in~\cite{CFLS}.

{\bf Proof of Theorem~\ref{main}.}
Let $\alpha=\alpha(4)$ be the constant from the Ramsey multiplicity problem above.  Let $\epsilon=10^{-9}$ and choose $c$ such that 
 every red/blue coloring of ${[N] \choose 2}$ for $N>c$ results in at least $(\alpha-\epsilon)N^4$ monochromatic copies of $K_4$. 
We will prove that
$$r(P_a^{2}, P_b^{2}) \le c(ab)^{9.7435}$$ 
for all $a,b \ge {2}$. This immediately gives the bounds we seek by letting $a=b=n$.

We will proceed by induction on $a+b$. If $a=2$ then the trivial upper bound is $b< c(2b)^{9.7435}$ and the same holds if $b=2$. For the induction step,  suppose we have a red/blue coloring of ${[N]\choose 2}$ with $N=c(ab)^{9.7435}>c$. By the choice of $c$, we obtain at least  $(\alpha-\epsilon) N^{4}$ monochromatic copies of $K_{4}$.

 Assume without loss of generality that half of these copies are  red. To each such copy with vertex set $x_0<x_1<x_2<x_3$ associate the middle pair of vertices $x_1, x_2$. By the pigeonhole principle there exists a set $Y=y_1<y_2$ that is the middle pair for at least $((\alpha-\epsilon)/2)N^{4}/{N-2 \choose 2}>(\alpha-2\epsilon)N^2=\beta N^2$ red copies of $K_{4}$. 
Let $L$ be the set of smallest vertices $y_0$ in these red $K_{4}$s and let $R$  be the set of largest vertices $y_{3}$ in these red $K_{4}$s. Note that  all edges between $Y$ 
and $L \cup R$ are red and that $|L||R|$ is the number of the red copies of $K_{4}$ that we are working with. 
If $|L|=\gamma N < \beta N$, then, since $L$ and $R$ are disjoint, $|R|\le (1-\gamma)N$ and so 
$$\beta N^2\le |L||R|\le \gamma(1-\gamma)N^2<\beta(1-\beta)N^2.$$
Consequently, $|L|\ge \beta N$ and similarly $|R|\ge \beta N$. By induction and $\log_2(1/\alpha)<9.7434$,
$$|L| \ge \beta N = \beta c (ab)^{9.7435}
> c (1/2)^{9.7435} (ab)^{9.7435}\ge r(P_{\lfloor a/2\rfloor}^{2}, P_b^2).$$

 Now we apply the induction hypothesis to find  either a red $P_{\lfloor a/2\rfloor }^{2}$
or a blue $P_b^{2}$ in $L$. If the latter occurs we are done so we get the former.  The same argument applies to $R$. 
Therefore, we obtain two red copies of $P_{\lfloor a/2\rfloor}^{2}$, one in $L$ and the other in $R$. Consider these two copies together with $Y$. Since the distance between vertices in $L$ and vertices in $R$ is at least $3$, and $Y$ is connected to all vertices in $L \cup R$ by red edges, this yields a red copy of $P_{2\lfloor a/2 \rfloor +2}^{2}$ which contains a red copy of $P_a^{2}$ because $2\lfloor a/2 \rfloor +2 \ge a$. \qed

\section{Proof of Theorem~\ref{main2}}

\begin{definition} \label{broomdef} The (ordered) broom $B^k_{a,m}$ is the $k$-graph with vertices $v_1 < v_2 < \cdots < v_a < w_1 <\cdots <w_m$
  such that $v_1, \ldots, v_a$ is a tight $k$-graph path and we also have all the edges $v_{a-k+2}\cdots v_aw_j$ for all $j \in [m]$.
  \end{definition}

We will omit the superscript $k$ in $B^k_{a,m}$ in all future usage as it will be obvious from the context. For example, in the proof below $k=3$. 
\bigskip
  
  \begin{theorem} \label{4n} $r_3(4, P_n) < 6nm$ where $m= r(P_n^2, P_n^2)$.
  \end{theorem} \proof  Recall that we are using the notation $P_n$ for a 3-uniform tight path and $P_n^2$ for the square of a 2-uniform (i.e. graph) ordered path. We will prove that every red/blue coloring $\chi$ of ${[N] \choose 3}$, where $N=6nm$ yields either a red $K_4^3$, or a blue $P_n$. Assume that there is no red $K_4^3$. We will show that there there is a blue $P_n$
     or a blue $B_{a, m}$ for all $2 \le a \le n$ within  the first $6am$ vertices. Since $P_n \subset B_{n,m}$ this will prove the result.  Let us show that there is a blue $P_n$ or a blue $B_{a, m}$ in $[6am]$ by induction on $a$. For the base case $a=2$, we seek a pair of vertices $v_1<v_2$ and at least $m$ vertices $w>v_2$ within $[12m]$ such that $v_1v_2w$ is blue. If we cannot find these $m$ vertices for any pair $v_1, v_2$, then the number of blue edges is at most ${12m \choose 2}(m-2) <  1/4{12m \choose 3}$ so the number of red edges is more than $(3/4){12m \choose 3}$ and a simple averaging argument shows that we would have a red $K_4^3$, contradiction. 
  
  Now for the induction step, assume that we have a blue copy of $B_{a-1, m}
 $ in $[6(a-1)m]$ and we wish to augment this to a blue copy of $B_{a, m}$ in $[6am]$. Suppose that the vertex set of the blue $B_{a-1, m}$ is
$$v_1 < v_2 < \cdots < v_{a-1} < w_1 <\cdots <w_m.$$
 Define the red/blue coloring $\phi$ of the complete graph  on $\{w_1, \ldots, w_m\}$ by $\phi(w_iw_j)=\chi(v_{a-1}w_iw_j)$.  By definition of $m$, we get  a copy $H$ of a monochromatic $P_n^2$ under $\phi$ with vertices $z_1<\cdots <z_n$.
 
  Suppose $H$ is  red under $\phi$. The four vertices $v_{a-1},z_i,z_{i+1}, z_{i+2}$ have three red edges, so $\chi(z_iz_{i+1}z_{i+2})$ is blue for all $i$. We conclude that 
 $z_1<\cdots <z_n$ is a blue $P_n$.  
 
 Next,  suppose $H$ is blue under $\phi$. Fix $i$ and consider the three vertices $z_i,z_{i+1}, z_{i+2}$.  If there are at least $m$  edges $z_iz_{i+1}y$ with $y \le 6am$ and $\chi(z_iz_{i+1}y)$ is blue, then we can use these edges to form a blue copy of $B_{a, m}$ in $[6am]$ with vertices
$$v_2<\cdots <v_{a-1} < z_i < z_{i+1} < Y$$
where $Y$ is the set of these $y$.
 So the number of such $y$ is at most $m$, and the same is true for the pairs $z_{i+1},z_{i+2}$
 and $z_i, z_{i+2}$. Since $6am-6(a-1)m=6m>3m$ there is a vertex $y$ such that $\chi(z_iz_{i+1}y)=\chi(z_iz_{i+2}y)=\chi(z_{i+1}z_{i+2}y)$ and these are all red.
 Therefore $\chi(z_iz_{i+1}z_{i+2})$ is blue. Since this argument applies for each $i$, we obtain a blue $P_n$ under $\chi$ with vertices $z_1<\cdots <z_n$. \qed
\bigskip

{\bf Proof of Theorem~\ref{main2}}. The case $k=3$ of Theorem~\ref{main2} follows immediately from Theorem~\ref{main}  and Theorem~\ref{4n}.  The lower bound for general $k$ follows from the lower bound for $r_k(P_{k+1},P_n)$ in \cite{MSW, MSh} (see also \cite{DLR}). The upper bound for general $k$ follows from the upper bound when $k=3$ (as a base case) and the standard pigeonhole argument for hypergraph Ramsey numbers due to Erd\H os and Rado (see the proof on pages 421-423 in~\cite{ER} for the original argument, or~\cite{GRS}, or Section 2 of~\cite{CFS}, or Section  2 of~\cite{MSEH}). 
Applying this argument from $k=3$ to $k=4$ raises the exponent of $n$ by a factor of 3 (from slightly less than $20.5$ to slightly less than $61.5$) and subsequent applications do not affect the exponent of $n$. \qed 
\bigskip

\section{Proof of Theorem~\ref{t=k}}

Let $f_k(n)$ be the minimum $N$ such that every red/blue coloring of 
${[N] \choose k}$ results in a blue $P_n$ or a set $S$ of $k+1$ vertices with at least $3$ red edges, one of which consists of the smallest $k$ vertices in $S$.
Let $H_k(3)$ denote the set of ordered $k$-graphs with three red edges as described above. We will abuse notation by saying a copy of $H_k(3)$ when we mean a copy of some $H \in H_k(3)$. 
\medskip

\begin{theorem}\label{generalkt=3}
$f_k(n) \le 2n^2$ for all $n >k\ge 3$.
\end{theorem}
\proof We will prove that every red/blue coloring $\chi$ of ${[N] \choose k}$, where $N=2n^2$ yields either a red $H_k(3)$  or a blue $P_n$. Assume that there is no red $H_k(3)$. We will show that there is a  blue $P_n$ or a blue  
 $B_{n-1, n}$. Since $P_n \subset B_{n-1, n}$ this will prove the result.  We will prove that there is a blue $B_{n-1, n}$ by showing that there is a blue $B_{a, n}$ in $[2an]$ for all $1 \le a \le n-1$ by induction on $a$.  The base case $a=1$ is trivial since $B_{1,n}$ is just a collection of $n+1$ vertices (with no edge since $k \ge 3$) and $2n>n+1$. 

  Now for the induction step, assume that we have a blue copy $B$ of $B_{a-1, n}$ in $[2(a-1)n]$ and we wish to augment this to a blue copy of $B_{a, n}$ in $[2an]$. Suppose that the vertex set of the blue $B_{a-1, n}$ is
$$v_1 < v_2 < \cdots < v_{a-1} < w_1 <\cdots <w_n.$$
Let $V=\{v_1, \ldots,v_{a-1}\}$ and $W=\{w_1, \ldots, w_n\}$.  For each $0 \le j \le k$, let $S_j=\{v_{a-j}, \ldots, v_{a-1}\}$ denote the $j$ largest vertices of $V$.   

{\bf Claim.} For every  $i \in [k]$ and $P \in {W \choose i}$, $\chi(S_{k-i} \cup P)$ is blue.

{\bf Proof of Claim.}
 Let us proceed by induction on $i$. The base case $i=1$ is trivial, due to the definition of $B$. Indeed, if $a-1 \ge k-1$, then all $k$-sets of the form $S_{k-1} \cup \{w\}$ for $w \in W$ are blue. If $a-1<k-1$ there is nothing to check. For the induction step, let $i \ge 2$ and suppose for contradiction that
 $\chi(S_{k-i} \cup P)$ is red for some  $P \in {W \choose i}$. 
 
  Let $P_1, P_2$ be two distinct $(i-1)$-sets contained in $P$. Note that $i \ge 2$ means that that ${i \choose i-1}\ge 2$ so such $P_1, P_2$ exist. Suppose that there are $n$ vertices $w_1'<w_2'<\cdots <w_n'\le 2an$ such that $w_1'>w_n$ and $\chi(S_{k-i} \cup P_1 \cup \{w_i'\})$ is blue.  Then $V \cup P_1 \cup \{w_1',\ldots, w_n'\}$
 is a blue $B_{a+i-2, n}$ which contains a blue $B_{a, n}$ in $[2an]$ as required. Indeed, it suffices to check that $\chi(S_{k-j} \cup Q_j)$ is blue for all $j \le i-1$, where $Q_j$ is the set of $j$ smallest vertices of $P_1$. But $Q_j \in {W \choose j}$ and $j<i$ so this is true by induction on $i$. We conclude that the number of such vertices $w_i'$ is at most $n-1$ and the same assertion holds for $P_1$ replaced by $P_2$. This gives at most $2n-2$ vertices between $w_n$ and $2an$. Since $w_n \le 2(a-1)n= 2an-2n$, there exists a $w$ such that $w_n<w\le 2an$ such that both $\chi(S_{k-i}\cup P_1 \cup \{w\})$ and
 $\chi(S_{k-i}\cup P_2 \cup \{w\})$ are red. This means that we have a red copy of $H_k(3)$ in $S_{k-i} \cup P \cup \{w\}$, contradiction.

  Now we simply apply the Claim with $i=k$ to conclude that all $k$-sets of ${W \choose k}$ are blue, and this is a blue clique with $n$ vertices which contains a blue $P_n$. \qed

\bigskip

{\bf Proof of Theorem~\ref{t=k}.} 
Let $f_k(k+1, t; P_n)$ be the minimum $N$ such that every red/blue coloring of 
${[N] \choose k}$ results in a blue $P_n$ or a set $S$ of $k+1$ vertices with at least $t$ red edges, one of which consists of the smallest $k$ vertices in $S$.
We observe that 
$f_{k}(k+1, t; P_n) < 2^{{f_{k-1}(k, t-1; P_{n-1}) \choose k}}$. Indeed, this follows from the standard pigeonhole argument for hypergraph Ramsey numbers due to Erd\H os and Rado (see the proof on pages 421-423 in~\cite{ER} for the original argument, or~\cite{GRS}, or  Section 2 of~\cite{CFS} or Section  2 of~\cite{MSEH}).
When applying this argument, we must note that the initial red edge in the $(k-1)$-graph gives rise to two red edges in the $k$-graph, one of which is again an initial edge. We apply this recurrence $t-3$ times until we have a $(k-t+3)$-graph.  As $k \ge t$, Theorem~\ref{generalkt=3} now applies to give  $f_{k-t+3}(k-t+4, 3; P_n)=f_{k-t+3}(n)\le 2n^2$ and this yields the result.  \qed

\section{Proof of Theorem~\ref{eh}}
Let us first prove the $t=3$ case of Theorem~\ref{eh} by improving the quadratic bound in Theorem~\ref{generalkt=3} to a linear bound.
 Let $F(3)$ be the collection of ordered $k$-graphs with $k+1$ vertices and at least three edges.  
  
\bigskip
  {\bf Proof of Theorem~\ref{eh} (a).}
  We are to show that $r_k(k+1, 3; P_n) \le 16n$ for all $n \ge k \ge 3$. 
   We will prove that every red/blue coloring $\chi$ of ${[N] \choose 3}$, where $N=16n$ yields either a red $H\in F(3)$ or a blue $P_n$.
Assume that there is no red $H \in F(3)$. We will show that there is a blue  $B_{n-1, 6}$.   As $P_n \subset B_{n-1,6}$ this will prove the result. We will prove that there is a blue $B_{n-1,6}$ by showing that
  there is a blue $B_{a, 6}$ for all $k-1 \le a < n$ within  the first $16(a+1)$ vertices by induction on $a$. For the base case $a=k-1$, we seek $k-1$ vertices $v_1<\cdots < v_{k-1}$ and at least $6$ vertices $w>v_{k-1}$ within $[16k]$ such that $v_1\cdots v_{k-1}w$ is blue. If we cannot find these $6$ vertices for any $(k-1)$-set $v_1<\cdots < v_{k-1}$, then the number of blue edges in $[16k]$ is at most $5{16k \choose k-1}$. A short calculation shows that this is less than $((k-1)/(k+1)){16k \choose k}$ since $k \ge 3$. Consequently, the number of red edges is more than 
$(2/(k+1)){16k \choose k}$ and an easy averaging argument then implies that   
 there is a  red $H \in F(3)$, contradiction.

  Now for the induction step, assume that we have a blue copy of $B_{a-1, 6}
 $ in $[16a]$ and we wish to augment this to a blue copy of $B_{a, 6}$ in $[16(a+1)]$. Suppose that the vertex set of the blue $B_{a-1, 6}$ is
$$v_1 < v_2 < \cdots < v_{a-1} < w_1 <\cdots <w_6.$$
For $q \in \{k-1, k-2, k-3\}$, let $S_q=\{v_{a-q}, \ldots, v_{a-1}\}$ be the greatest $q$ vertices among the $v_i$.
 Define the red/blue coloring $\phi$ of the complete graph  on $\{w_1, \ldots, w_6\}$ by $\phi(w_iw_j)=\chi(S_{k-2}w_iw_j)$.  Since $r(3,3)\le 6$, we obtain a monochromatic triangle $T$ under $\phi$. If $T$ is red, then 
  $T \cup S_{k-2}$ yields a red member of $F(3)$ which is a contradiction, so $T$ is blue. Assume for simplicity that $T=w_1<w_2<w_3$. 
Fix $1 \le i <j \le 3$.    If there are at least $6$  edges $S_{k-3}w_iw_jy$ with $w_3<y \le 16(a+1)$ and $\chi(S_{k-3}w_iw_jy)$ is blue, then we can use these edges to form a blue copy of $B_{a, m}$ in $[16(a+1)]$ with vertices
$$v_2<\cdots <v_{a-1} < w_i < w_j < Y$$
where $Y$ is the set of these $y$.
 So the number of such $y$ is at most $5$, and the same is true for all three pairs  $\{i,j\}$. Since $16(a+1)-16a=16>15$ there is a vertex $y$ such that $\chi(S_{k-3}w_1w_2y)=\chi(S_{k-3}w_1w_{3}y)=\chi(S_{k-3}w_2w_3y)$ and these are all red. This gives a red member of $F(3)$, contradiction. \qed
\bigskip

 {\bf Proof of Theorem~\ref{eh} (b).}
For the lower bound, consider the ordering of $2n-3$ points
$$v_1, v_1' < v_2, v_2'< \cdots < v_{n-2}, v_{n-2}' < v_{n-1}.$$
For all $i<j$, color the triples  $v_iv_i' v_j$ red. Color all other triples blue. No blue path $P_q$ can contain both $v_i$ and $v_i'$ unless they appear at the end of $P_q$. If no such pair is in $P_q$ then clearly $q \le n-1$.
If $v_i, v_i'$ are both in $P_q$, then they lie at the end of the path, and $q \le (i-1)+2=i+1 \le n-1$. Hence there is no blue $P_n$. If we have four points with two red edges, then two of these triples must be $v_iv_i'a$ and $v_iv_i'b$ for some $v_i<a,b$. But then there cannot be any other red edge among these points.

For the upper bound, we will prove the following stronger statement by induction on $n$:
$$r_3(4,3; B_{n-1,2})\le 3n-3.$$

Since $P_n \subset B_{n-1,2}$, this will prove the upper bound. 
The base case $n=3$ is true due to the following simple argument. Suppose we have a red/blue coloring of ${[6] \choose 3}$. If there are two edges of the form $12i$ that are blue then we have a blue $B_{2,2}$ so there is at most one such edge. This means that there are $2<i<j<k\le 6$ such that  $12x$ is red for all $x \in S=\{i,j,k\}$. If any edge of the form $1xy$ or $2xy$ is red  for $x,y \in S$, then $\{1,2,x,y\}$ has three red edges, so all such edges are blue. This means that $1ij$ and $1ik$ are both blue, giving a blue $B_{2,2}$. 

 Now for the induction step, assume that $N=3n-3$ and we have a 
red/blue coloring $\chi$ of the triples of $[N]$ with no four points containing three red edges.  By induction, we obtain a blue $B_{n-2,2}$ in $[N-3]$. Say the vertices of this $B_{n-2,2}$ are 
$$v_1< \cdots < v_{n-2} < v_{n-1}< v_{n-1}'.$$ Let $S=\{N-2, N-1, N\}$. If there are two blue edges of the form $v_{n-2}v_{n-1}x$ and $v_{n-2}v_{n-1}y$
 for $x,y \in S$, then we have obtained a blue copy of $B_{n-1,2}$ with $x=v_n$ and $y=v_n'$. So at most one of these  edges is blue. The same applies to $v_{n-2}v_{n-1}'x$ and $v_{n-2}v_{n-1}'y$. Since $|S|=3$, there is $w \in S$ such that both $v_{n-2}v_{n-1}w$ and $v_{n-2}v_{n-1}'w$ are red.
  Now consider the four vertices 
$v_{n-2}, v_{n-1}, v_{n-1}', w$. We have identified two red edges among them so $v_{n-2}v_{n-1}v_{n-1}'$ and $v_{n-1}v_{n-1}'w$ are blue. If there exists $z \in S-\{w\}$ such that $\chi(v_{n-1}v_{n-1}'z)$ is blue, then we obtain a copy of $B_{n-1,2}$ as follows:
$$v_2< \cdots < v_{n-2} < v_{n-1}< v_{n-1}'< w <z \qquad {\rm or} \qquad
v_2< \cdots < v_{n-2} < v_{n-1}< v_{n-1}'< z<w.$$
Therefore $\chi(v_{n-1}v_{n-1}'a)=\chi(v_{n-1}v_{n-1}'b)$ and both are red, where $S=\{a,b,w\}$ with $a<b$. If $\chi(v_{n-2}v_{n-1}x)$ is blue for some $x\in \{a,b\}$, then we obtain the blue $B_{n-1,2}$
$$v_1<\cdots < v_{n-2} < v_{n-1} < v_{n-1}'< x.
$$
Hence $\chi(v_{n-2}v_{n-1}x)$ is red for both  $x\in \{a,b\}$. The four vertices $v_{n-2}, v_{n-1}, v_{n-1}', x$ contain the two red edges $v_{n-2}v_{n-1}x$ and
$v_{n-1}v_{n-1}'x$ so  $\chi(v_{n-2}v_{n-1}'x)$ is blue for $x \in \{a,b\}$. This gives us
$$v_1<\cdots < v_{n-2} < v_{n-1}' <  a < b$$
which is a blue $B_{n-1,2}$, and the proof of $r_3(4,3; B_{n-1,2})\le 3n-3$ is complete.

Observe that $r_3(4,3; P_n)\le r_3(4,3; B_{n-1,2})-1$ by  taking an optimal construction for $r_3(4,3; P_n)$, adding a new largest vertex $v$ and  coloring all triples containing $v$ with blue. \qed
\bigskip

 {\bf Proof of Theorem~\ref{eh} (c).}
We are to show that $r_4(5, 4; P_n) \ge 2^{n-2}+1$. Let us proceed by induction on $n$. The case $n=2$ is trivial, so assume we have a construction for $n-1$ that uses the vertices $[2^{n-3}]$. To obtain the construction for $n$, take a copy of the construction for $n-1$ among the vertices $\{2^{n-3}+1, \ldots, 2^{n-2}\}$. It remains to color 4-sets that intersect both halves of $[2^{n-2}]$. 
We color all the 4-sets that have exactly two points in each half red and all other $4$-sets blue. Let us first argue that no 5 points contain 4 red edges. If all 5 points lie in one half then we are done by induction.  If the distribution of points is $4+1$ then we again have at most one red edge and if the distribution is $3+2$ then we have exactly three red edges. The other cases are of course symmetric. Next we argue that there is no blue $P_n$. Such a blue $P_n$ cannot have two points in both halves as the $2+2$ edges are red so all but one point must lie in one half. This gives a blue $P_{n-1}$ in one half which cannot exist by induction.  \qed
\bigskip

\noindent \textbf{Acknowledgment.}  
I am grateful to Andrew Suk for many helpful discussions on this topic, to David Conlon for informing me about the results in~\cite{BCKK, CFLS}, and to both referees who read the paper very carefully and provided many comments that improved the presentation.

\end{document}